\documentclass[12pt]{article} 
\title{\bf Shot Noise with Random Parameters}
\author{Jean-Fran\c{c}ois Chamayou\thanks{Laboratoire de
Statistique et Probabilit\'es, Universit\'e Paul Sabatier, 31062 Toulouse, France. e-mail: \texttt{jfchamay@math.univ-toulouse.fr}}}

\date{\today}
\usepackage[cp850]{inputenc}
\usepackage[T1]{fontenc}
\usepackage[english]{babel}
\usepackage{graphicx}

\def\reel{\hbox{I\hskip -2pt R}}

\def\nne{\hbox{I\hskip -2pt N }}

\begin{document}
\maketitle

{\bf Abstract}:
 
   We give the description of the following model: 
  $$  U_{n}=X_{n}(Y_{n}+U_{n-1})$$ for $n>1$
in the case where the  $X_{n}$ are i.d.d. random variables with probability density:
$$ A x^{A-1} , x \in [0,1] ,$$
    $A$ is also a random variable distributed according to a Gamma law.
    The $ Y_{n}$ are or deterministic and equal to $1$ or independent Gamma random variables.
    We use this model to compute the shot noise with random parameters.\\


AMS CLASSIFICATION :
Primary 60F

 \section{\bf Introduction:}
In many field and particularly in Physics the following process occurs:
\begin{eqnarray*} Z ( t ) = \Sigma_{0 \le T_i \le t} Y_{i} \exp{(- B ( t - T_i ))} \end{eqnarray*}
where the $ (T_i)_{i \in \nne^{*} } $ are the arrival dates of a Poisson process with rate $ \Lambda $  and a decay parameter $B$, $Y_{i}$ a deterministic amplitudes equal to $1$ or independent random amplitudes with  Laplace transform $ w_{Y}(s)$, i. e.  $Z(t)$ represents the shot noise, the stationary limit of this process
exists and is the same as the distribution limit of the sequence of random variables $ ( U_n)_{n \in \nne^{*} } $
 defined by 
\begin{eqnarray*} U_{1}=X_{1}Y_{1},U_{n+1}=X_{n+1}(Y_{n+1}+U_{n}) , n > 0 \end{eqnarray*} 
The $ X_n $ being powers in $ \frac{B}{\Lambda} $ of the independent variables $ V_n $ uniform on$ [0,1]$,
the $Y_{n}$ are determistic $=1$ or independent Gamma varibles.
see Chamayou and Dunau [6,7,9] .
$ f_{U} (u) $ the probability density of the stationary case with deterministic amplitude is given by the following first order delay differential equation, see Chamayou and Dunau [10]:
\begin{eqnarray*} u f^{'}_{U} (u) + f_{U} (u) = - \frac{\Lambda}{B} ( f_{U} (u) - f_{U} (u-1)) , u>0 \\
 f_{U} (u) = 0 , u < 0 \end{eqnarray*}
  $ g_{U,Y} (s)$ the Laplace transforms are given by, see Chamayou and Schorr [5]et [8]:
\begin{eqnarray*} g_{U,Y} (s) = \exp{( - \frac{\Lambda}{B}  \int^{s}_{0} \frac{ (1 - h_{Y}(\xi))}{\xi} d\xi )} \end{eqnarray*} 
and the Laplace transform $w_{Y}(s)=\exp{(-s)}$ for the deterministic amplitude.
we deduce by succesive integrations on the intervals $ [n, n+1], n \geq 0 $
\begin{eqnarray*} f_{U} (u) = \frac{ \exp{(- \gamma  \frac{\Lambda}{B} )} u^{ \frac{\Lambda}{B} -1 } } { \Gamma(  \frac{\Lambda}{B} ) } , 0 \le u \le 1 \\
 f_{U} (u) =  \frac{ \exp{(- \gamma \frac{\Lambda}{B}) } u^{ \frac{\Lambda}{B} -1 } } { \Gamma(  \frac{\Lambda}{B} ) } ( 1 -
 \int_{1}^{u} \frac{ (\eta -1)^{\frac{\Lambda}{B}-1} } { \eta ^{ \frac{\Lambda}{B} } } d\eta ), 1 \le u \le 2 \end{eqnarray*}
where $ \gamma $ is the Euler constant and $ \Gamma $ the  Gamma function.  $ F_{U} (u) $ the distribution function is given by:
\begin{eqnarray*} F_{U} (u) = \frac{( \exp{- \gamma  \frac{\Lambda}{B} )} u^{ \frac{\Lambda}{B} } } { \Gamma(  \frac{\Lambda}{B}+1 ) } , 0 \leq u \leq 1 \\
F_{U} (u) =  \frac{ \exp{(- \gamma  \frac{\Lambda}{B} )} u^{ \frac{\Lambda}{B}} } { \Gamma(  \frac{\Lambda}{B}+ 1 ) } ( 1 -
 (u-1)^{\frac{\Lambda}{B} +1 } \int_{1}^{u} (\frac{ (\eta -1)} { \eta })^{ \frac{\Lambda}{B} }  d\eta ), \\
 1 \leq u \leq 2 \end{eqnarray*}

\section{\bf Random Parameters:}

\subsection{\bf General Case:}
For the case where $ B $ and/or $ \Lambda $ are random, denote by $ A = \frac{\Lambda}{B} $ the random variable with probability density $ h_A (a) $ independent of the $ (V_n)_{n \in \nne^{*} } $ we get for the Laplace transform:
\begin{eqnarray*}
g_{U,A,Y} (s) =\int_{0}^{\infty} \exp{( - a \int^{s}_{0} \frac{ (1 - w_{Y}(\xi) )}{\xi} d\xi )} h_A (a) da \end{eqnarray*} 

\subsection{\bf Gamma Case:}
For the case where  $ A $ is distributed according to a Gamma law with parameter $ \alpha $ , we get the Laplace transform:
\begin{eqnarray*} 
g_{U,A,Y} (s) =\int_{0}^{\infty} \exp{( - a (1+ \int^{s}_{0} \frac{ (1 - w_{Y}(\xi))}{\xi} d\xi))} \frac{a^{\alpha -1}}{\Gamma(\alpha)} da \\
 g_{U,A,Y} (s) = \frac{1}{(1+ \int^{s}_{0} \frac{ (1 - w_{Y}(\xi))}{\xi} d\xi )^\alpha}
 \end{eqnarray*}
The behaviour for $ 0< u \ll 1$ for the deterministic amplitude is given from:
$$ \lim_{s \rightarrow \infty }  g_{U,A} (s) = \frac{1}{(1+ \gamma + \ln {s} )^\alpha} $$
since see Abramovitz [1]
$$
\int_{s}^{\infty} \frac{\exp{(-\xi)}}{\xi} d\xi = - \gamma - \ln{s} + \int^{s}_{0} \frac{ (1 - \exp{( - \xi)} )}{\xi} d\xi $$
i.e. : 
$$ \frac{1}{\ln^{\alpha}{(\exp{(1+\gamma)} s)}}$$
is for the Laplace transform of the distribution function:
$$ \frac{1}{s \ln^{\alpha}{(\exp{(1+\gamma)} s)}} $$
is the Laplace transform of the Volterra $ \mu $ function:
$$ \mu ( \exp{(-(1+\gamma))} u , \alpha - 1) ,$$ 
see Apelblat [2,3,4],  which equals the Volterra $ \nu $ function: 
$$ \nu( \exp{(-(1+\gamma))} u ) $$
for the exponential case: $ \alpha = 1 $ we have taken as computational example , where:
\begin{eqnarray*} \nu( z )  = \int_{0}^{\infty} \frac{ z^{\xi} }{ \Gamma(\xi+1) } d\xi \end{eqnarray*}

The asymptotic behaviour $ 1 \ll u $ is given from:
$$ \lim_{s \rightarrow 0 }  g_{U,A} (s) = \frac{1}{(1+ s )^\alpha} $$
which gives rise to a gamma law with parameter $\alpha $ .

\subsection{\bf Exponential Case $ \alpha = 1 $ :}

In the exponential case it is easy from:
$$ \lim_{s \rightarrow 0 }  g_{U,A} (s) = \frac{1}{(1+ s - \frac{s^2}{4})} $$
and from the Laplace transform inversion:
$$\aleph(u)= \frac{(\exp{(-(2\sqrt{2}-2)u)}-\exp{(-(2\sqrt{2}+2)u)})}{\sqrt{2}} $$

we could have $ N > 2$  terms in the expansion of 
$ \int^{s}_{0} \frac{ (1 - \exp{( - \xi)} )}{\xi} d\xi $
for the price of the computation of $N$ roots  to obtain the sum of  $ N $ exponentials.

The Laplace transform of the distribution function: 
$$ \frac{1}{s(1+ s )^\alpha} $$

gives rise to an incomplete gamma function $ \gamma(\alpha,s) $, see [1].

 In the same way for the distribution function from:
$$ \lim_{s \rightarrow 0 }  g_{U,A} (s) = \frac{1}{s(1+ s - \frac{s^2}{4})} $$
we get:
$$ \Xi(u)= 1-\frac{(\frac{\exp{(-(2\sqrt{2}-2)u)}}{(2\sqrt{2}-2)}-\frac{\exp{(-(2\sqrt{2}+2)u)}}{(2\sqrt{2}+2)})}{\sqrt{2}} $$

\section{\bf Exponential Parameters Example ( $\alpha = 1$ ):}

We get for the probability density:
\begin{eqnarray*} f_{U} (u) = \int_{0}^{\infty} \frac{ \exp{(- (\gamma +1) a )} u^{ a -1 } } { \Gamma(  a ) } da
=\\
\frac{1}{u} \phi( \exp{(-(1+\gamma))} u ) , 0 \leq u \leq 1 \end{eqnarray*}
where the Fransen-Wrigge function $ \phi$  [12] reduces to: 
$$\phi(z) = \int_{0}^{\infty}  \frac{ k( z^{a} ) }{ \Gamma(a) } da $$
where:
$$ z =  \exp{(-(1+\gamma))} u , k(z) = z $$
We will use for the numerical calculations the given weights and abscissae, numerical calculations consist in an integration of the Fransen and Wrigge function on the following interval:
\begin{eqnarray*} f_{U} (u) = 
\frac{1}{u} \phi( \exp{(-(1+\gamma))} u ) - \phi( \exp{(-(1+\gamma))}(u -1)) \\
+ \frac{1}{u} \int_{1}^{u}  
\phi( \exp{(-(1+\gamma))}u(1-\frac{1}{\xi})) d\xi ,
 1<u<2 
\end{eqnarray*}
see  figure $1$ representing the simulation of the process density of probability and the previously described numerical computation, the difference is due to the fact that the density tends to $ \infty$ for the variable $u$ tending to $0$. The approximation of the density is given by $ \aleph(u) $ for $ u > 2 $. \\

The distribution function equals:
\begin{eqnarray*} F_{U} (u)= \nu( \exp{(-(1+\gamma))} u ) , 0 \le u \le 1 \end{eqnarray*}
in that domain the asymptotical representation given by Wyman and Wong [13] used by  Apelblat has the adequate accuracy for a comparaison with the simulation:
\begin{eqnarray*} \nu (\xi) \approx \frac{-1}{\ln{\xi}}( 1 + \frac{ \psi(1) }{ \ln{\xi} } + \frac{ (\psi^{2}(1)-\psi\prime(1))}{\ln^2{\xi}} + ...)
\end{eqnarray*}
where the functions $ \psi^{(n)} , n \ge 0 $ are the logarithmic derivatives of the Gamma function, the numerical computation can be performed using the coefficients $ a_j , j= 1,..., 61 $ given by Fransen et Wrigge [11] :
\begin{eqnarray*} \nu(\xi) \approx \sum_{j=0}^{N} \frac{a_j j!}{(-\ln{\xi})^{j+1}} \end{eqnarray*}
the expansion is stopped at a $ N $ derived from the argument $ \xi $ and
\begin{eqnarray*} F_{U} (u) = 
\nu( \exp{(-(1+\gamma))} u ) - (u-1) \nu ( \exp{(-(1+\gamma))}(u -1)) \\
 +  \int_{1}^{u}  
\nu ( \exp{(-(1+\gamma))}u(1-\frac{1}{\xi})) d\xi,
 1<u<2 
 \end{eqnarray*}
For the arguments of the function $ \nu $ superior to  $ \exp{(-(1+\gamma))}$, the computation is performed by integration as in Apelblat. The approximation of the distribution function is the following: 
$  F_{U} (2) + \Xi(2) - \Xi(u) $ pour $ u > 2 $.
See figure $2$ representing the simulation of the process distribution and the result of the above numerical computation.

\section{\bf Random Amplitudes:}
\subsubsection{\bf Gamma Case:}

For a random amplitude $Y$ Gamma distributed, the Laplace transforms is given by:$ w_{Y}(s) =\frac{1}{(1+s)^{\beta}}$ 

\subsubsection{\bf Exponential case $\beta = 1$:}
The Laplace transform for $ \beta = 1 $ is given by: 
$$ g_{U,A,Y}(s) =  \frac{1}{\ln ^{\alpha} {(\exp{(1)}(1+s))} } $$
From the Laplace transform inversion in the exponential case for $A , \alpha = 1 $ the probability density is easily deduced using the Fransen and Wrigge function:
\begin{eqnarray*}
 f_{U} (u) = \exp{(-u)} \int_{0}^{\infty} \frac{ \exp{(-  a )} u^{ a -1 } } { \Gamma(  a ) } da
=\\
\frac{\exp{(-u)}}{u} \phi( \exp{(-1)} u ) , u > 0 
\end{eqnarray*}
See  figures $3$ and $4$ representing the simulation of the process probability density and  distribution function compared to the result of the numerical computation of the above formula, its integration requires the use of the approximate value for
 $ u \rightarrow 0 $:
$$ \exp{(-u)}\nu( \exp{(-1)} u) .$$
\subsubsection{\bf Case $\beta= 1/2$:}
From the density for $\beta=1/2$ using the parabolic cylinder function $D$ given in [6]:
\begin{eqnarray*}
 f_{U,A} (u) = \sqrt{\frac{2}{\pi}} \exp{(-(\frac{u}{2}))} 2^{3 A} A u^{A-1} D_{-(1+2 A)} (\sqrt{2 u}) , u >0
\end{eqnarray*}
We deduce from the integral representation (see [1]) changing the order of integrations the probability density using the Fransen et Wrigge function, this avoids an integration on the index:
\begin{eqnarray*}
 f_{U} (u) = \frac{\exp{(-u)}}{4u}\sqrt{\frac{2}{\pi}} \int_0^{\infty} \exp{-(\sqrt{2u}\xi+\frac{\xi^2}{2})}
\phi(\exp{(-\frac{1}{2})} 2 \sqrt{2 u} \xi) d\xi , u>0
\end{eqnarray*}
See  figure $5$ representing the simulation of the process probability density compared to the result of the numerical computation of the above formula

\subsubsection{\bf Case $\beta= 2$:}
The density for $\beta=2$ is given using the modified Bessel function of first kind $I$ given in [6]:
\begin{eqnarray*}
 f_{U,A} (u) = \exp{(-(A+u))}  (\frac{u}{A})^{\frac{A-1}{2}} I_{A-1} (2\sqrt{A u}) , u >0
\end{eqnarray*}
using the series expansion of this function, one can deduce (see [1]) 
\begin{eqnarray*}
 f_{U} (u) = \exp{(-u-2)} \sum_{k=0}^{\infty} \exp{(2 k)} \mu (\exp{(-2)}u , k , k-1) , u > 0
 \end{eqnarray*}
where $ \mu (z,b,a) $ is one of the Volterra's functions, see Apelblat [2,3,4]

\subsection{\bf Laplace distributed Amplitudes:}
For both signs random amplitude the Fourier transform is used instead of Laplace transform, in the case of a Laplace distribution:
$$ \tilde{w}_{Y}(s) =\frac{1}{(1+s^2)^{\beta}}$$ 

\subsubsection{\bf $\beta=1$  Case:}
Then the Fourier transform for $ \beta = 1 $ is:
$$ \tilde{g}_{U,A,Y}(s) =  \frac{1}{\ln ^{\alpha} {(\exp{(1)}(1+s^2))} } $$
From the density for $\beta=1$ given in [5,6]using the modified Bessel function of second kind $K$  ( MacDonald function)
using its integral representation see [1],  we deduce the density in terms of Fransen and Wrigge function in the exponential 
case for $A , \alpha = 1 $:
\begin{eqnarray*}
 f_{U} (u) = \frac{2}{\sqrt{\pi}|u|} \int_0^{\infty} \exp{(-(\xi+\frac{u^2}{4 \xi}))}
\phi(\frac{\exp{(-2)}}{\xi}(\frac{u}{2})^2) \frac{ d\xi}{\sqrt{\xi}}  , u \in \reel
\end{eqnarray*}
The distribution function is deduced from formal inversion of the Fourier transform for positive values due to symetry:
\begin{eqnarray*}
 F_{U} (u) = \frac{2}{\pi} \int_0^{\infty}
 \frac{ \sin{\xi} }{\xi ( 1+\frac{1}{2} \ln{(1+ (\frac{\xi}{u})^2)})} d\xi , u>0
\end{eqnarray*}
See  figures $6$ and $7$ representing the simulation of the process probability density and distribution function compared to the result of the numerical computation of the above formulae
(the functions are computed on the absolute value of the random variable due to the symetry).

\subsubsection{\bf $\beta=2$  Case:}
The case $\beta=2$ gives for $u \in \reel$ a difference of 2 non-central $\chi^2 $  with $A$   degrees of freedom and a non-centrality parameter $A$, according to [8]:
\begin{eqnarray*}
 f_{U,A} (u) = 
 \frac{\exp{(-(\frac{A}{2}))}}{\sqrt{\pi}}( \frac{|u|}{2})^{\frac{(A-1)}{2}}\sum_{n=0}^{\infty}
 \frac{(A |u|)^n}{\Gamma(n+1) \Gamma(\frac{A}{2}+n) 2^{2n}}
  K_{\frac{(A-1)}{2}+n}(|u|) ,
 u \in \reel
 \end{eqnarray*}
then we have:
 \begin{eqnarray*}
 f_{U} (u) = \frac{2}{\sqrt{\pi}|u|}\sum_{n=0}^{\infty} (\frac{u}{2})^{2n} \int_0^{\infty} \exp{(-(\xi+\frac{u^2}{4 \xi}))}
\mu(\frac{\exp{(-3)}}{\xi}(\frac{u}{2})^2,n,n-1)\frac{ d\xi}{\sqrt{\xi}} , u \in \reel
\end{eqnarray*}
The distribution function is deduced from formal inversion of the Fourier transform for positive values due to symetry:
\begin{eqnarray*}
 F_{U} (u) = \frac{2}{\pi} \int_0^{\infty}
 \frac{ \sin{\xi} }{\xi ( 1+\frac{1}{2} (\ln{(1+ (\frac{\xi}{u})^2)}+\frac{\xi^2}{(\xi^2 + u^2)}))} d\xi , u>0
\end{eqnarray*}
See figure $8$ representing the simulation of the process distribution function compared to the result of the numerical computation of the above formula
(the distribution function  is computed on the absolute value of the random variable due to the symetry).

\subsubsection{\bf $ \beta = \frac{1}{2} $ Case:}
Then the Fourier transform for $ \beta = \frac{1}{2} $ is: 
$$ \tilde{g}_{U,A,Y}(s) =  \frac{1}{\ln ^{\alpha} {(\exp{(1)}\frac{(1+\sqrt{(1+s^2)})}{2})} } $$
The distribution function is deduced from formal inversion of the Fourier transform for positive values due to symetry:
\begin{eqnarray*}
 F_{U} (u) = \frac{2}{\pi} \int_0^{\infty}
 \frac{ \sin{\xi} }{\xi ( 1+\ln{\frac{(1+\sqrt{(1+ (\frac{\xi}{u})^2)})}{2}})} d\xi , u>0
\end{eqnarray*}
See figure $9$ representing the simulation of the process distribution function compared to the result of the numerical computation of the above formula
(the distribution function is computed on the absolute value of the random variable due to the symetry) .
\newpage
\begin{center}
\begin{picture}(400,200)(0,0) 

\put(200,0){\vector(1,0){100}}
\put(10,0){\line(0,1){200}}

\thicklines
\put(0,0){\line(1,0){200}}

\put(50,50){\oval(80,80)[bl]}
\put(100,90){\oval(100,100)[bl]}
\put(190,110){\oval(180,180)[bl]}
\put(260,160){\oval(140,140)[bl]}
\put(285,180){\oval(50,40)[bl]}

\thinlines
\put(145,190){rate$ =\Lambda$}
\put(145,170){decay$ =B$}
\put(1,10){...$t_{0}$}
\put(50,10){$t_{1}$}
\put(100,10){$t_{2}$...}
\put(160,10){...$t_{n-k}$...}
\put(260,10){...$t_{n}$}
\put(270,110){$Z(t)$}

\put(8,0){$\uparrow$}
\put(48,0){$\uparrow$}
\put(98,0){$\uparrow$}
\put(258,0){$\uparrow$}

\put(10,0){\line(0,1){200}}
\put(50,0){\line(0,1){200}}
\put(100,0){\line(0,1){200}}
\put(190,0){\line(0,1){200}}
\put(260,0){\line(0,1){200}}

\put(280,0){$t$}
\put(20,150){Scheme 1: Shot Noise Process}

\end{picture}
\end{center}

\newpage

\begin{figure}
\centerline{\includegraphics[width=15cm]{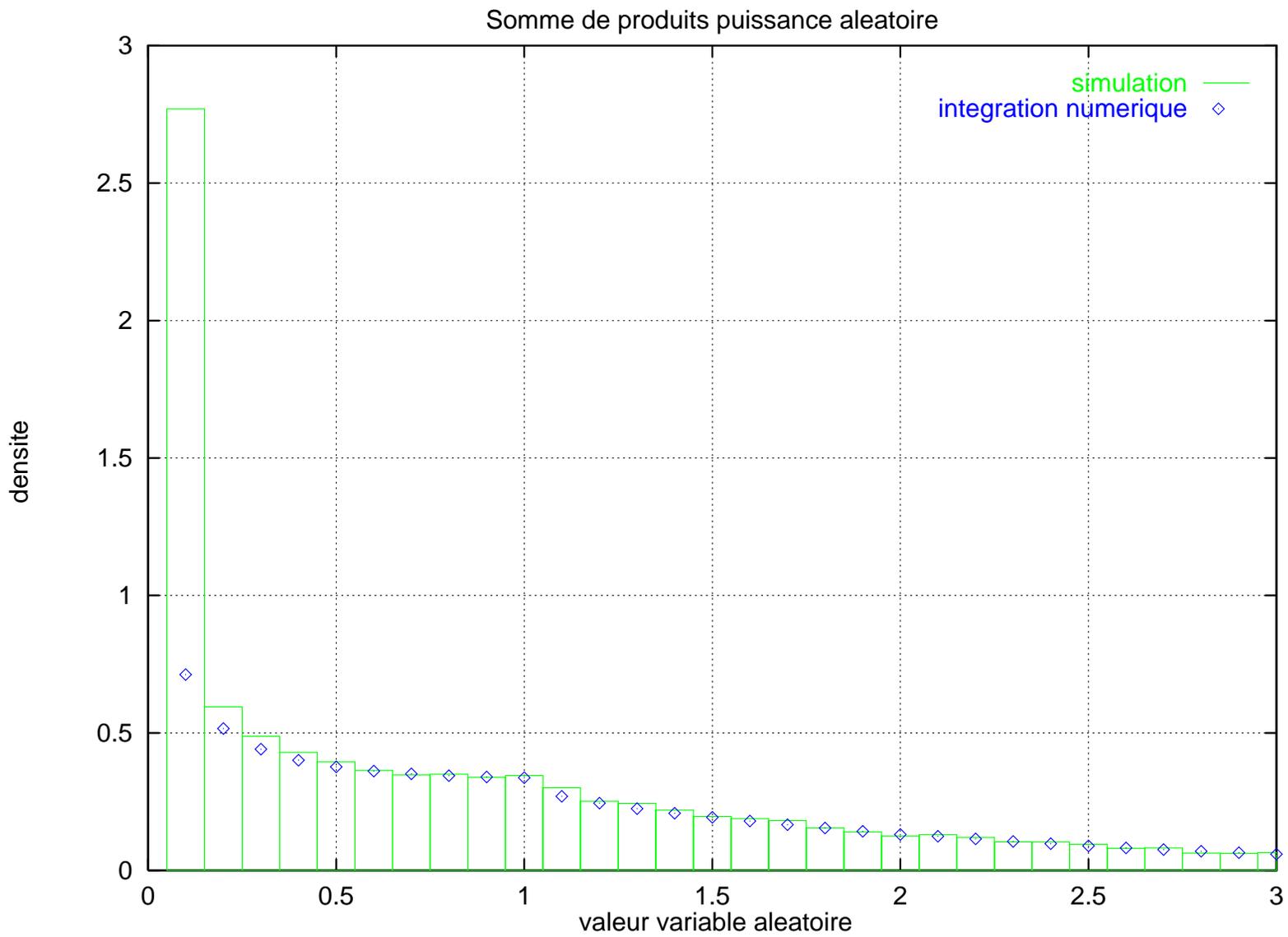}}
\caption{probability density (sum of random powers products)}
\end{figure}
\newpage 

\begin{figure}
\centerline{\includegraphics[width=15cm]{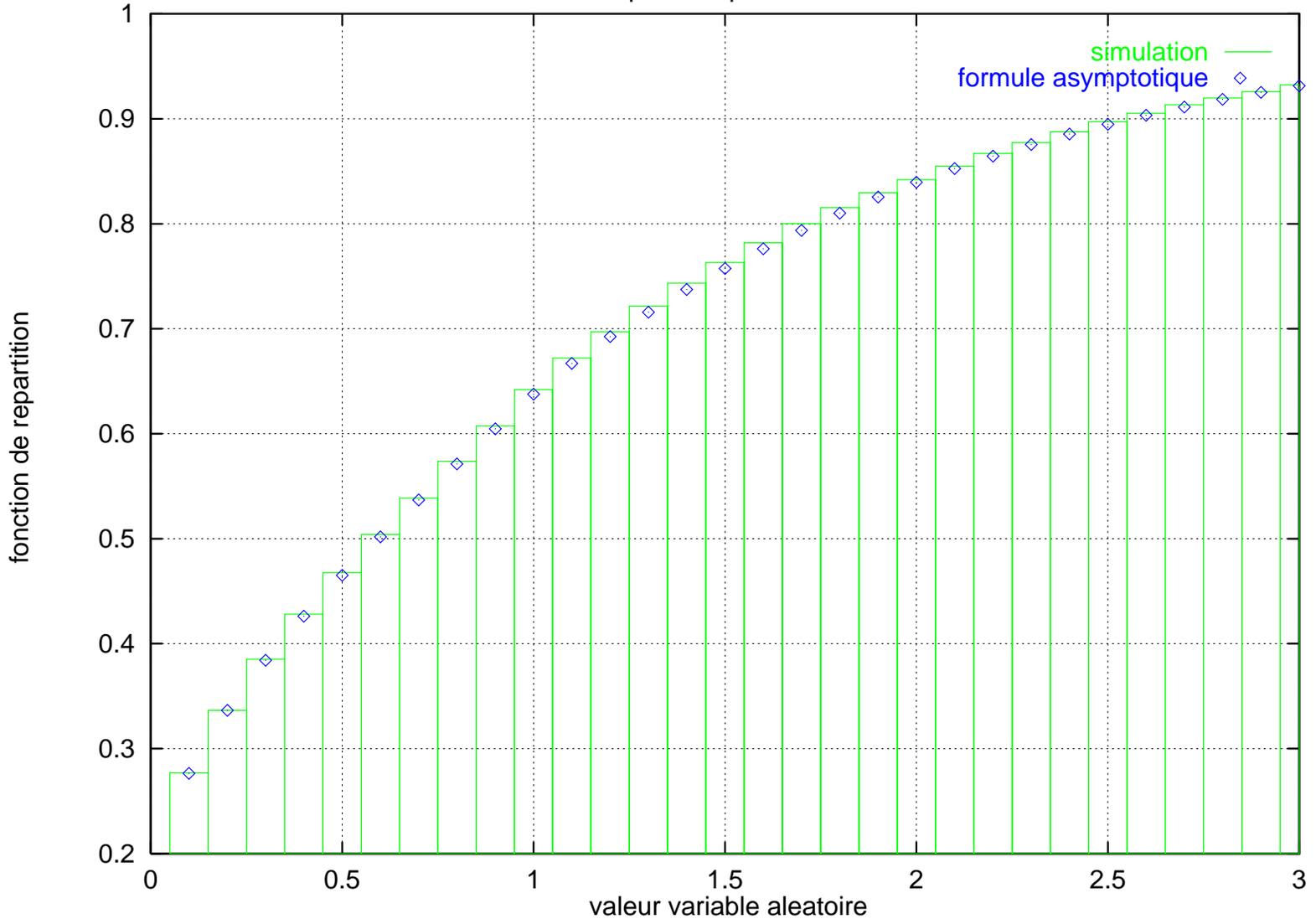}}
\caption{distribution function (sum of random powers products)}
\end{figure}

\begin{figure}
\centerline{\includegraphics[width=15cm]{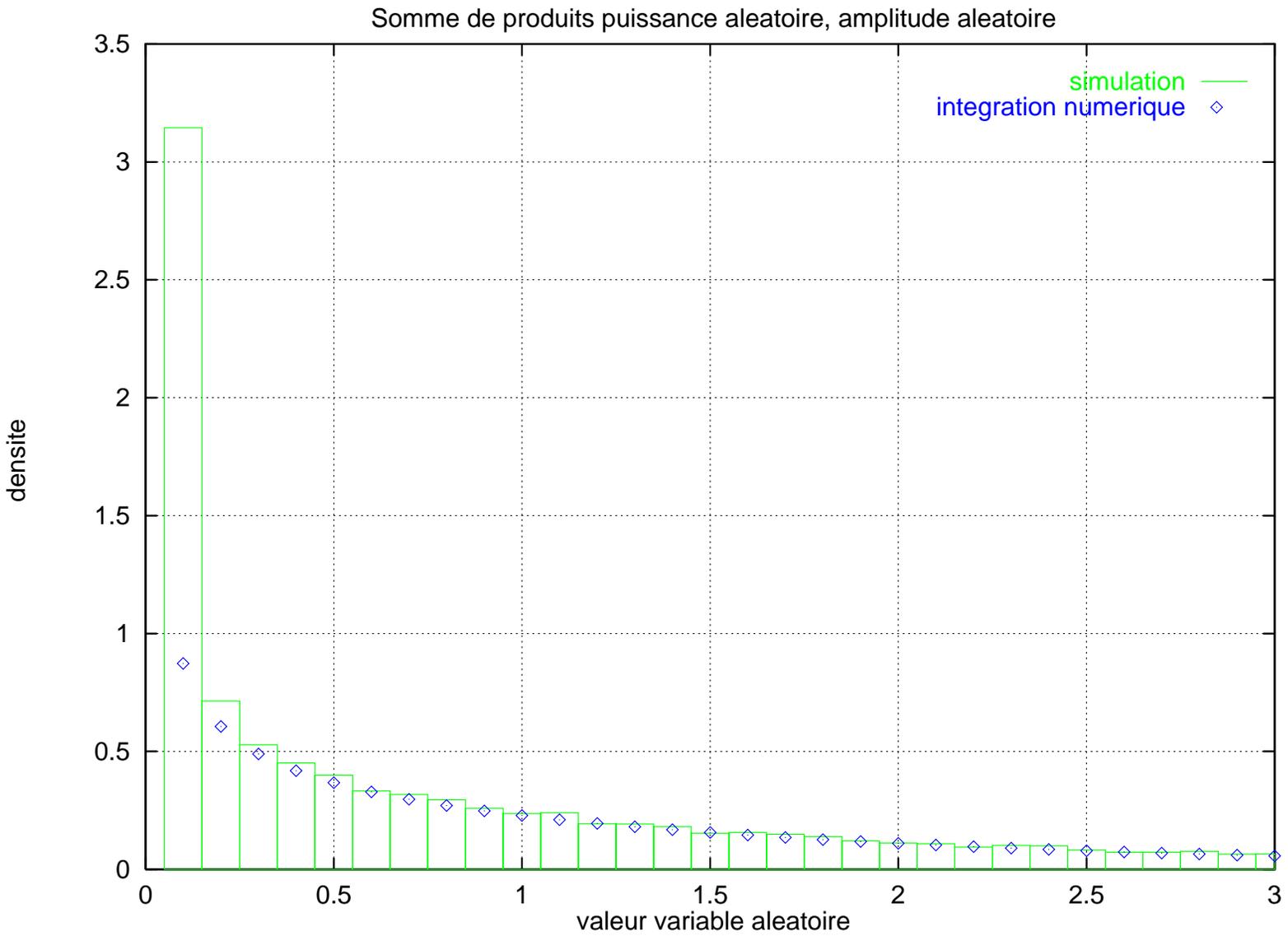}}
\caption{ probability density (sum of random powers products with exponential random amplitude $1$)}
\end{figure}
\newpage 

\begin{figure}
\centerline{\includegraphics[width=15cm]{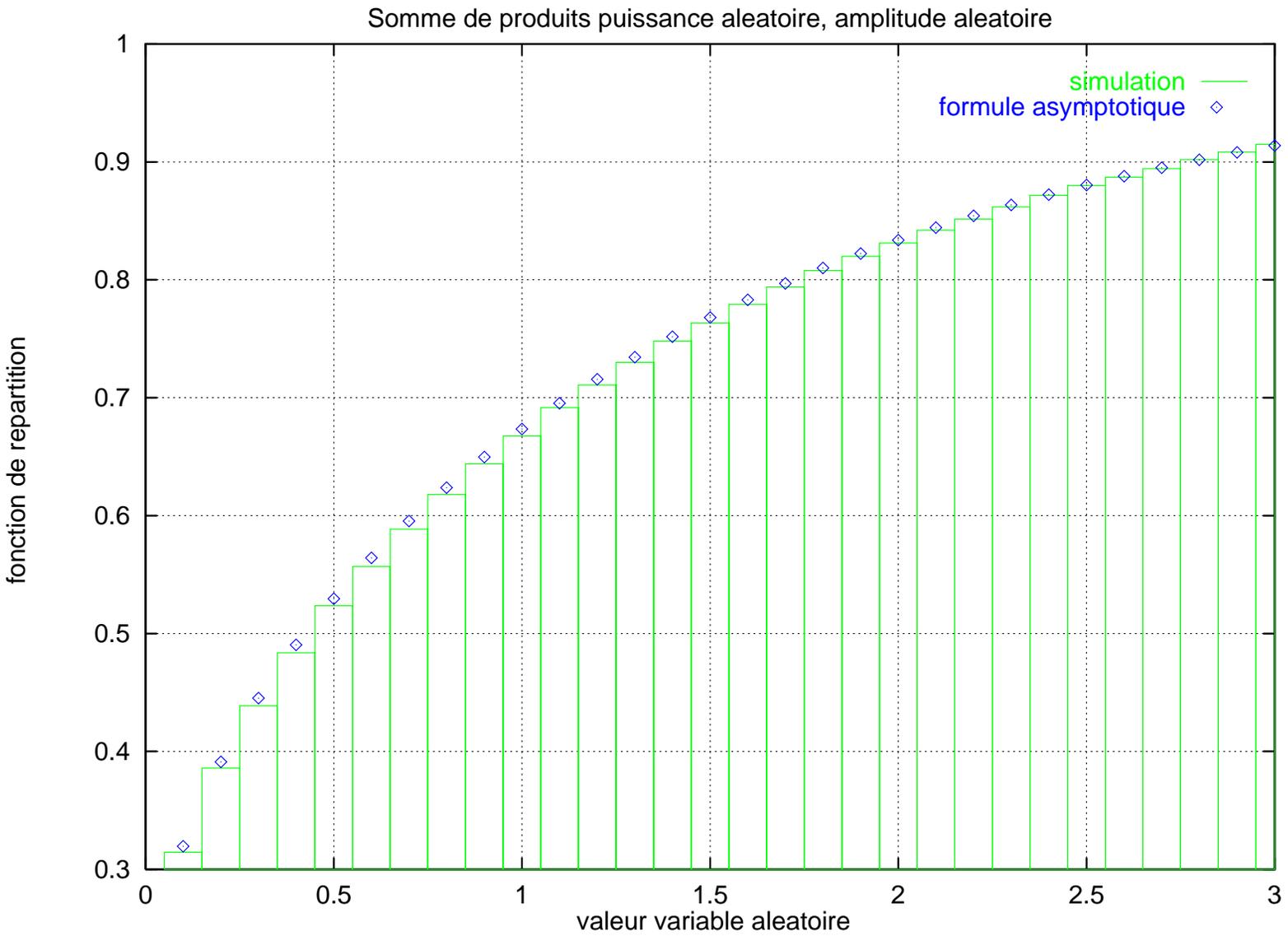}}
\caption{distribution function (sum of random powers products with exponential random amplitude $1$)}
\end{figure}

\begin{figure}
\centerline{\includegraphics[width=15cm]{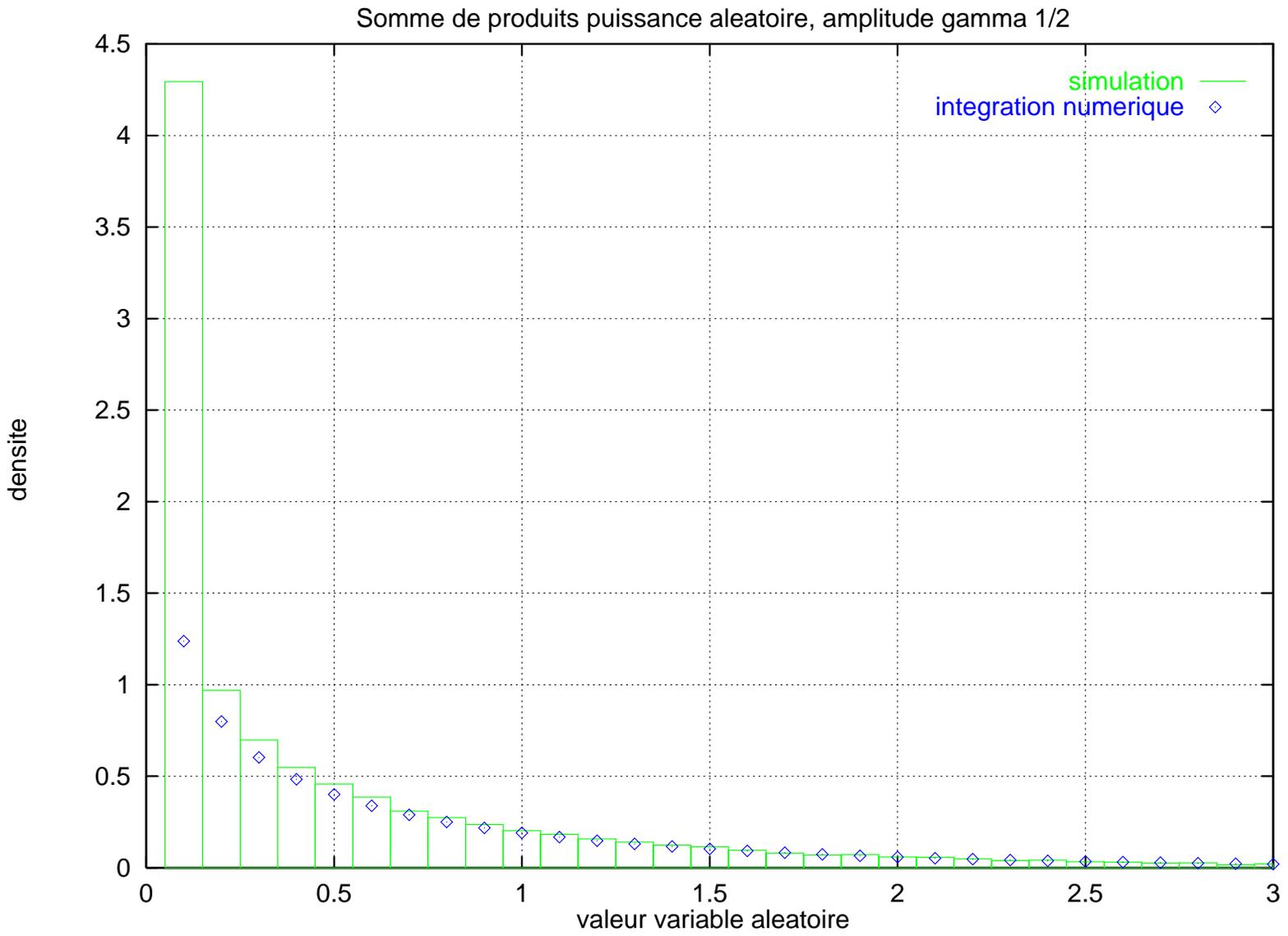}}
\caption{ probability density ( Gamma distributed random amplitude $\beta= 1/2$}
\end{figure}
\newpage 
 \begin{figure}
 \begin{center}
\includegraphics[width=15cm]{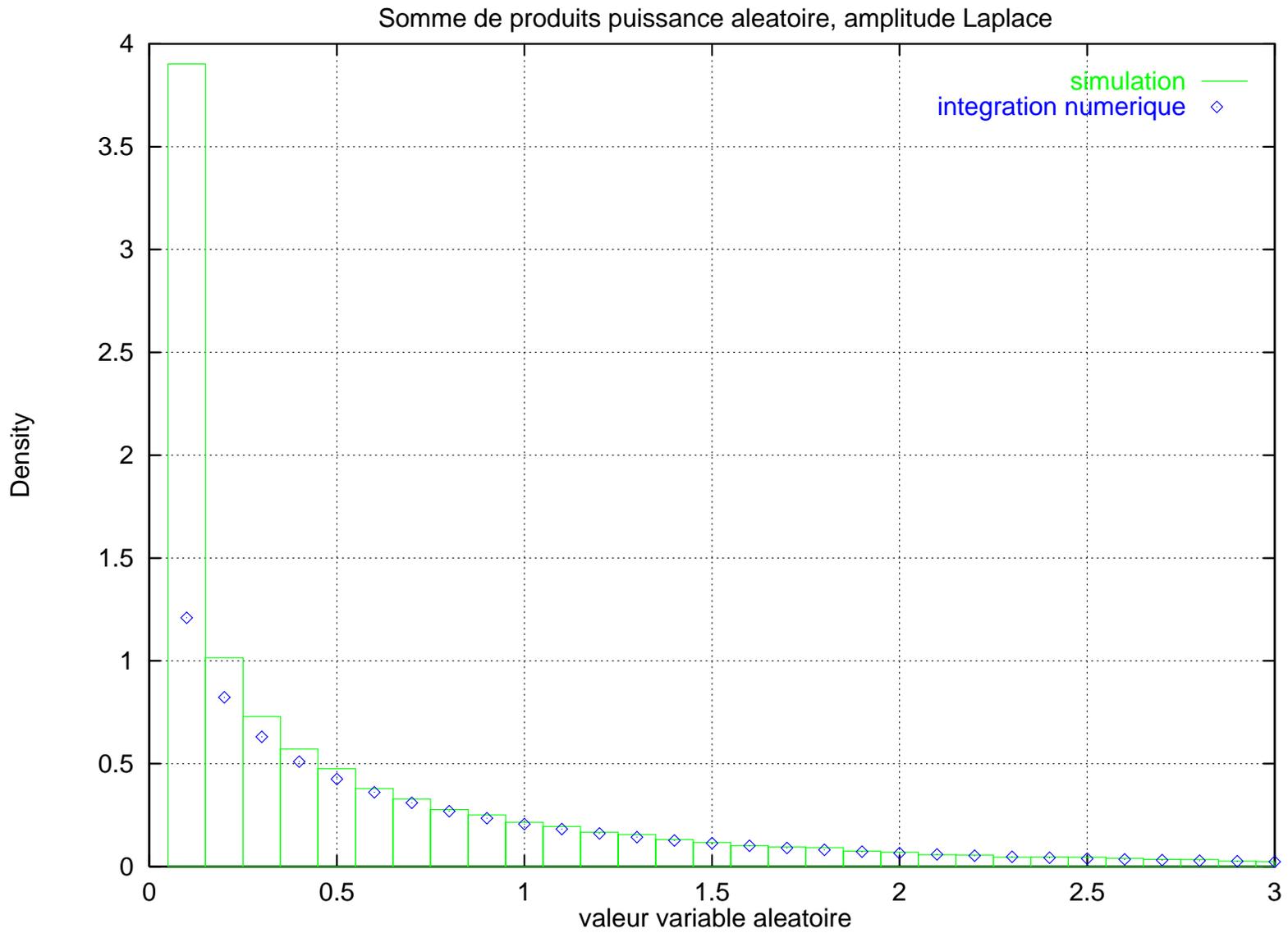}
\caption{probability density, Laplace distributed Amplitude $1$}
 \end{center}
\end{figure}
\newpage
\begin{figure}
\centerline{\includegraphics[width=15cm]{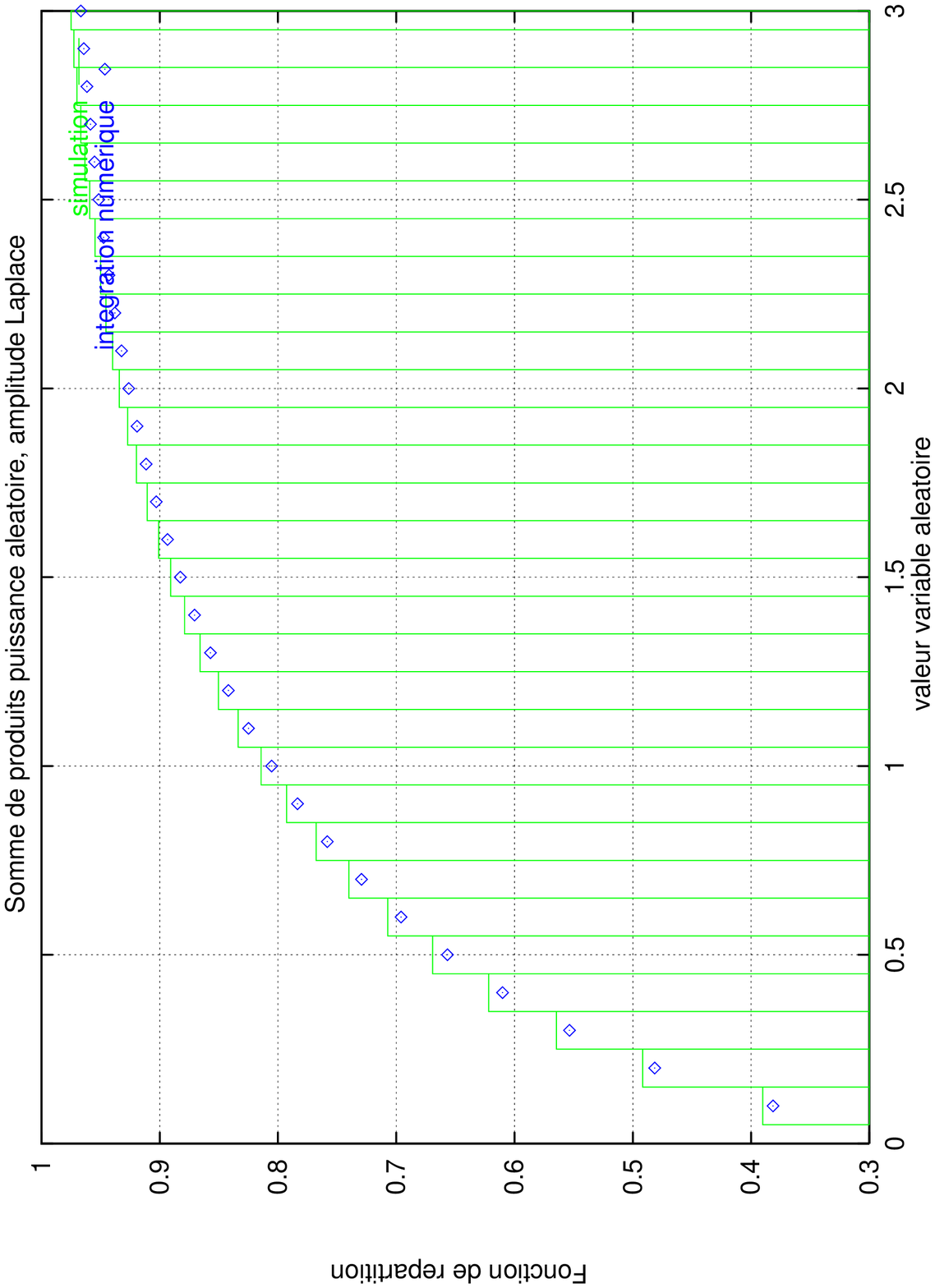}}
\caption{distribution function, Laplace distributed Amplitude $1$}
\end{figure}
\newpage
\begin{figure}
\centerline{\includegraphics[width=15cm]{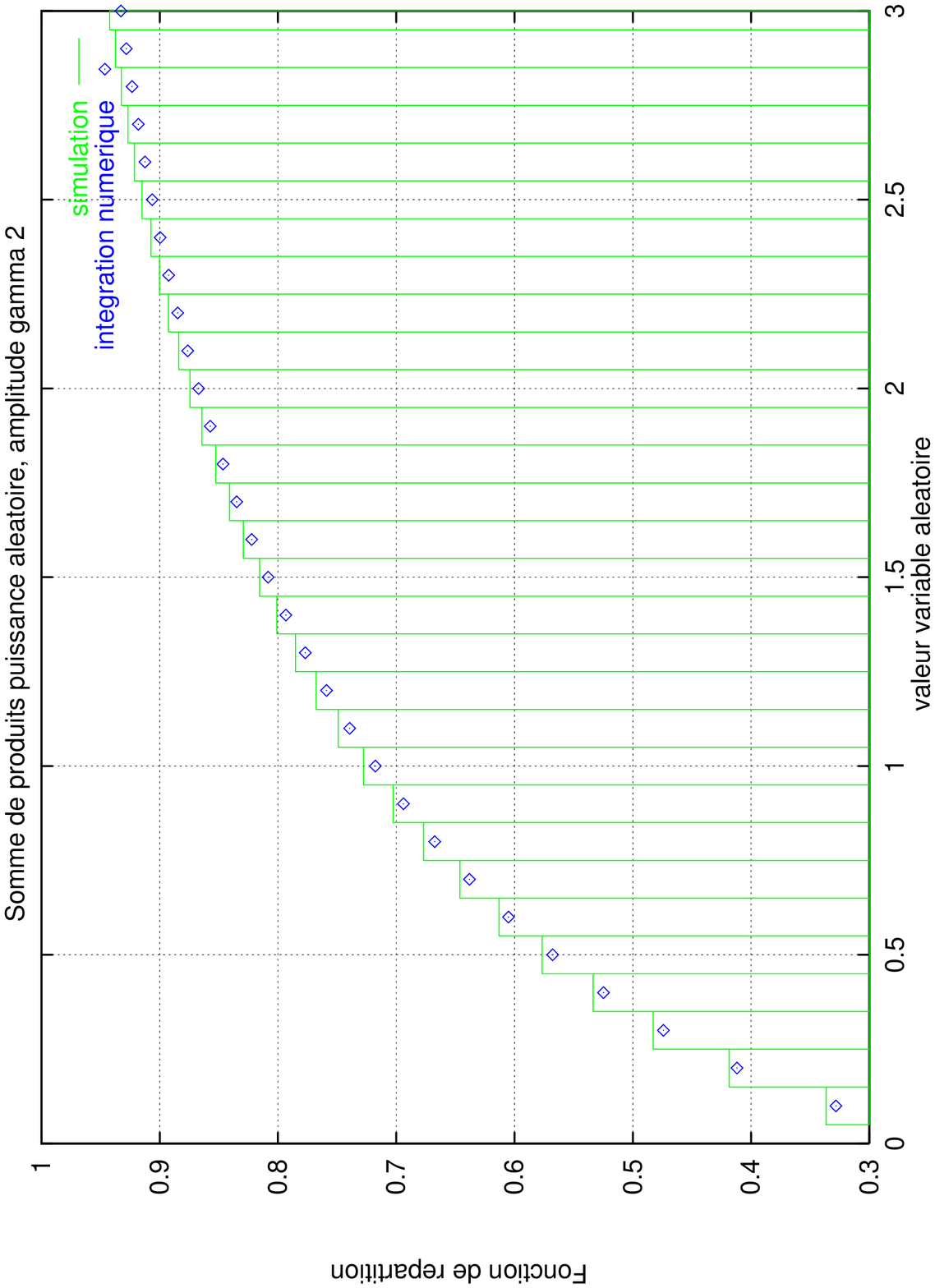}}
\caption{distribution function, Laplace distributed Amplitude $2$}
\end{figure}
\newpage
\begin{figure}
\centerline{\includegraphics[width=15cm]{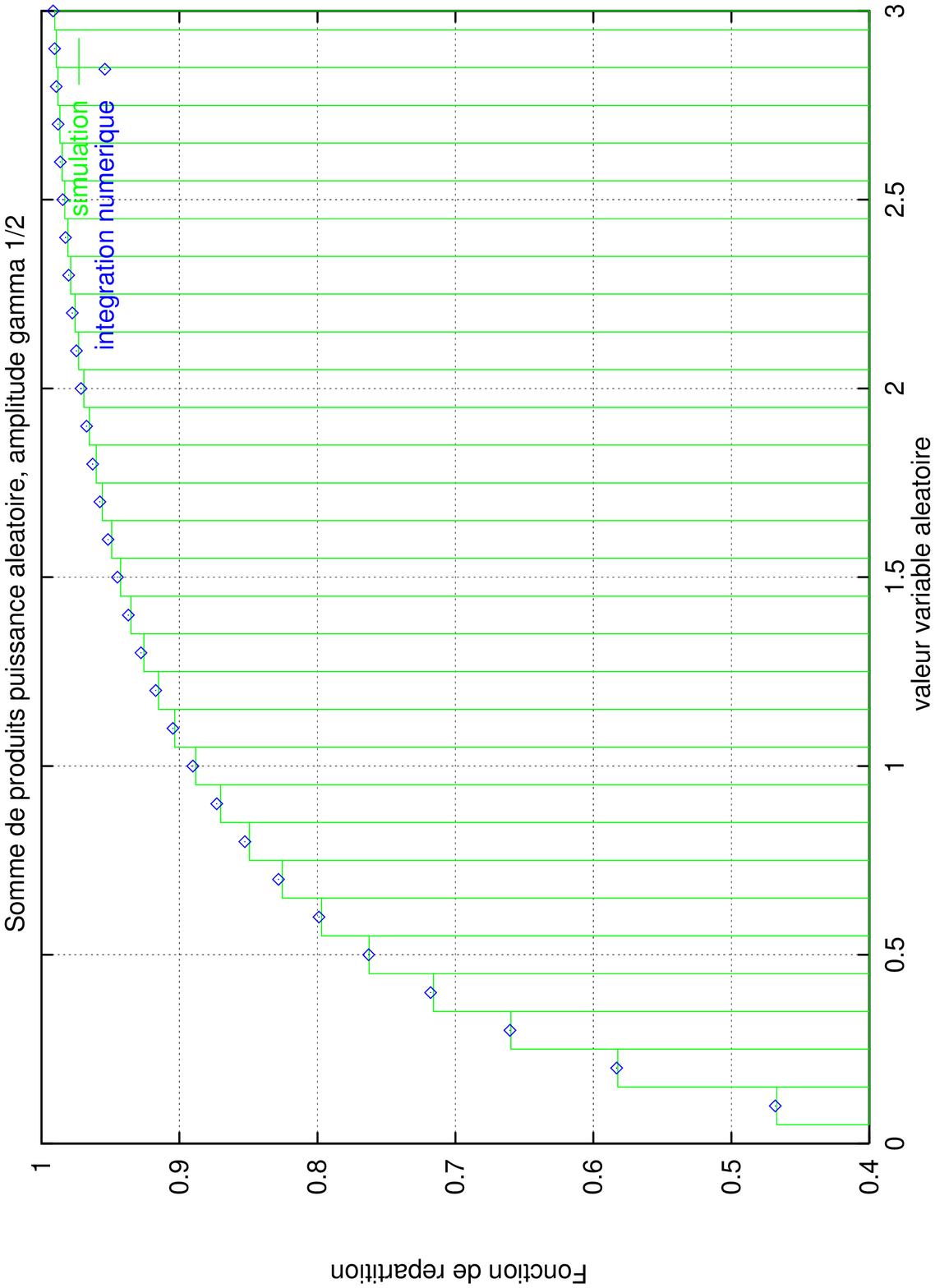}}
\caption{distribution function, Laplace distributed Amplitude $1/2$}
\end{figure}

\end{document}